\section{REPRESENTATIVE NUMERICAL RESULTS}
\label{Sec:S4_Numerical_results}
In the following section we present results from three numerical examples 
involving flow and transport in porous domains with heterogeneous material
properties. For all examples, the mesh was constructed from linear triangular 
elements. Identical meshes were used for both the RT0 and VMS elements. 
A direct solver (Amesos umfpack \cite{amesos:para06}) was used to solve 
the resulting system of equations in serial on a desktop machine, using an 
Intel Xeon CPU (2.93 GHz). For all examples, a constant species diffusivity, 
$D = 0.01$ was used. Also, the Backward Euler method was used for time 
integration with a time step, $\Delta t = 0.01 \; \mathrm{s}$, and the residual 
tolerance was set to $1 \times10^{-9}$.

\subsection{Multilayer problem}
For the multilayer problem (which is a modification of an example given in \cite{Wan_Thesis, MasudDG}), the unit square domain $\Omega = (0,1) \times 
(0,1)$, shown in Figure \ref{fig:layers_Geo}, is divided into five layers of thickness 0.2 units, each with a unique ratio, $k/\mu$. Within each layer, the ratio is constant, but changes as a step function at the layer boundaries. The flow is driven by a pressure difference prescribed on the left and right boundary of the domain. The properties used for this problem are defined in Table \ref{table:MultilayerProblemParams}. The computational mesh used for the multilayer problem is shown in Figure \ref{fig:layers_Geo}.

The exact solution for the pressure is a linear drop across the domain. The exact solution for the velocity in each layer is equal to the ratio, $k/\mu$. Figure \ref{fig:Layers_Ux_vs_y} shows the $x$-velocity along the centerline of the domain taken at $x = 0.5$. Note from the figure that the RT0 element captures the exact solution for the velocity in each layer and the jump in velocity between layers. The VMS element is not able to accurately capture the jump in velocity between layers. This is due to the well-known property that a continuous element will enforce continuity of the velocity in both the normal and tangential directions, although discontinuity in the tangential direction is physically realizable. Ironically (for a \emph{stabilized} method), in addition to smoothing the velocity in the vicinity of permeability jumps, the VMS element shows oscillations in the velocity within each layer. It can be shown that with increasing disparity between layer properties, these oscillations increase in magnitude. 

Figures \ref{fig:Layers_CvsY}, \ref{fig:Layers_C_contour}, and \ref{fig:Layers_CvsTime} show the results of using the computed velocity field as the advection velocity for the transport problem defined by equation \eqref{Eqn:AD_Equilibrium}. In this case, the problem represents transport of a species though a stratified system of permeabilities. Along the left edge of the domain, the concentration is held constant in time at 1.0. As time progresses, the species advects and diffuses though the domain at different rates, given the local flow velocity with each layer. Even though the VMS element shows significant overshoots and undershoots for the velocity in the vicinity of the permeability jumps, the two elements show close correlation in terms of concentration of the species in time. Figure \ref{fig:Layers_CvsY} shows that the VMS element exhibits latency in tracking the species front in time, but Figure \ref{fig:Layers_CvsTime} suggests that the overall concentration, or the integral of $c$ over the domain is preserved. There is only a marginal difference between the two elements in predicting the time at which a steady state solution is reached. Also, both elements capture a smooth variation of concentration along the front, even thought the VMS velocity field is not smooth.

Figure \ref{fig:Layers_Mass_Balance} shows the degree to which the VMS element violates local mass balance for problems of this nature. (A similar figure is not presented for the RT0 element because it is locally conservative to machine precision.) The error was calculated as $\mathrm{div}[\boldsymbol{v}]$ over each element. Notice that the error is concentrated along the inflow and outflow boundaries. This is due to the well-known implications of imposing the pressure boundary conditions in weak fashion.
\begin{remark}
Due to the well-known instabilities associated with a high ratio of advection to diffusion in the transport model \cite{Masud2,TurnerGFEM}, some diffusion was necessary to include for numerical reasons, even though the primary quantities of interest are related to pure advection. Also a constant scalar diffusivity was chosen to avoid instabilities associated with violations of the discrete maximum-minimum principle \cite{Liska_Shashkov_CiCP_2008_v3_p852,Nakshatrala_Valocchi_JCP_2009_v228_p6726,
Nagarajan_Nakshatrala_IJNMF_2010}. 
\end{remark}
\begin{table}[htb!]
\caption{Multilayer problem: parameters}
\centering
\begin{tabular}{l l}
\hline
Parameter & Value\\
\hline
$\rho$ & 996.1 kg/$\mathrm{m}^2$\\
$k_1$ & 1$\times10^{-13}$ $\mathrm{m}^2$\\
$k_2$ & 5$\times10^{-13}$ $\mathrm{m}^2$\\
$k_3$ & 0.5$\times10^{-13}$ $\mathrm{m}^2$\\
$k_4$ & 3$\times10^{-13}$ $\mathrm{m}^2$\\
$k_5$ & 8$\times10^{-13}$ $\mathrm{m}^2$\\
$\mu$ & 1$\times10^{-8}$ kg/m-s\\
$\boldsymbol{g}$ & (0,0,0)\\
$p_o$ & 1$\times10^{5}$ Pa\\
\hline
Nodes & 441\\
Elements & 800\\
\hline
\end{tabular}
\label{table:MultilayerProblemParams}
\end{table}
\subsection{Cylinder inclusion problem}
The cylinder inclusion problem can be conceptualized as a contaminant initially contained in a highly permeable porous media, encapsulated in a lower permeability media, that is dispersed by a cross-flow in the domain generated by a pressure drop. This problem is an adaptation from the problem presented in \cite{Gupta}. The primary quantity of interest in this problem is the degree to which the contaminant infiltrates the domain and the time scale of this process. The geometry of this problem is shown in Figure \ref{fig:Inclusion_Geo} along with the computational mesh. The permeability of the infinite domain, and other parameters is given Table \ref{table:CylinderInclusionParams}.

Figure \ref{fig:Vel_x} shows the velocity magnitude for both elements. For the RT0 element, a smooth transition occurs from the highly permeable region to the region of lower permeability. On the other hand, the VMS element shows oscillations along the boundary of the two regions. Again, this is due to the continuity of the solution in the direction transverse to the interface for the VMS element. Figures \ref{fig:Pressure_vs_x} and \ref{fig:Vel_vs_x} show the pressure and velocity along different portions of the domain. The overshoot of the velocity magnitude for the VMS element is also evident in these figures.

Figure \ref{fig:Inclusion_C} shows the concentration of the species in time as it deviates from the initial conditions. Both the RT0 and VMS elements again show tremendous similarity in the contours of concentration with respect to time. The local mass balance error for the VMS formulation is shown in Figure \ref{fig:Inclusion_Mass_Balance}. Unlike in the multilayer problem, for the cylinder inclusion, the error is concentrated along the boundary of the inclusion, in particular where the interface is oriented tangentially to the flow direction. To more fully understand if these features are related to the magnitude of the jump in permeability, several cases were run in which the permeability of the inclusion was varied from 1$\times10^{-11}$ to 1$\times10^{-3}$ $\mathrm{m}^2$. Figure \ref{fig:Inclusion_Mass_Balance_Ratio} shows the error in the concentration for both elements as the jump increases. As expected, the RT0 element performs well regardless of the jump in permeability, but even the VMS shows only a marginal error (less than 2\%) even for the worst case scenario.
\begin{table}[htb!]
\caption{Cylinder inclusion problem: parameters}
\centering
\begin{tabular}{l l}
\hline
Parameter & Value\\
\hline
$\rho$ & 996.1 kg/$\mathrm{m}^2$\\
$k_1$ & 1$\times10^{-8}$ $\mathrm{m}^2$\\
$k_2$ & varies 1$\times10^{-11}$ to 1$\times10^{-3}$ $\mathrm{m}^2$\\
$\mu$ & 1.13$\times10^{-3}$ kg/m-s\\
$\boldsymbol{g}$ & (0,0,0)\\
$p_o$ & 2.1721$\times10^{5}$ Pa\\
\hline
Nodes & 489\\
Elements & 912\\
\hline
\end{tabular}
\label{table:CylinderInclusionParams}
\end{table}

\subsection{Leaky well problem}
The leaky well problem represents the transport of a species as it is injected into an aquifer that contains an abandoned well that allows some of the species to escape. The primary quantity of interest is the rate at which the species escapes from the aquifer through the leak. This problem has been studied by a number of authors in evaluating subsurface modeling codes (see \cite{Ebigbo,Nordbotten,Nakshatrala_Turner_IJES_2010}) and is therefore important to study from the perspective of element choice, since the element choice could have a significant effect on the results. The domain is shown in Figure \ref{fig:LeakyWellDomain} along with the computational mesh. 

The steady state pressure and velocity contours are shown in Figures \ref{fig:LeakyWellVMags}  and \ref{fig:LeakyWellP}. Given the velocity field computed by either element, the transport problem is solved by prescribing a Dirichlet boundary condition of $c = 1.0$ at the inflow. Thus, the injection rate is $\int cv_y \; \mathrm{dA}_i$ and the leak rate is $\int cv_y \; \mathrm{dA}_o$, where $\mathrm{A}_i$ is the area of the injection port and $\mathrm{A}_o$ is the area of the leaky well, both with diameter 0.3 m. For the leak rate, $c$ is measured at the opening of the leaky well into Aquifer A. Figures \ref{fig:RT0LeakyWellContours} and \ref{fig:VMSLeakyWellContours} show similar results for the overall contours of concentration for results obtained using the Darcy flow model given in equation \eqref{Eqn:FlowTransport_constant}.

Figure \ref{fig:LeakRate} shows the leak rate for various values of $\beta$ from 0.0 to 1$\times 10^{-9}$ using the Barus model given in equation \eqref{Eqn:FlowTransport_Barus}. Figure \ref{fig:LeakRate} shows that the dependence of the viscosity on the pressure has a strong influence on the predicted leak rate. This trend is also apparent in Figure \ref{fig:ContoursCompare}, which shows contours of the concentration at time $t = 0.12$ days for various values of $\beta$. Notice the greater degree of penetration for the low values of $\beta$ vs. high. The pressure dependent viscosity leads to both a lower leak rate and a longer time to steady-state conditions. Figure \ref{fig:LeakRate} reveals that the leak rate is substantially under-predicted by the VMS element. Although both elements predict that the leak rate will reach a steady-state around $t = 2$ days for the case of $\beta = 0.0$, the VMS element shows a considerable decrease in the steady-state leak rate. Whereas for the previous two problems, 
the RT0 and VMS element performed remarkably similar, in this case the VMS element shows considerable variation. The results suggest that the VMS element is not an appropriate choice for this type of problem since it not only shows a significant accuracy error, but also an under-prediction of the quantity of interest which, from a design standpoint, is less conservative. 

\begin{table}[htb!]
\caption{Leaky well problem: parameters}
\centering
\begin{tabular}{l l}
\hline
Parameter & Value\\
\hline
$\rho$ &  479 kg/$\mathrm{m}^2$\\
$k_1$ & 1$\times10^{-12}$ $\mathrm{m}^2$\\
$k_2$ & 1$\times10^{-14}$ $\mathrm{m}^2$\\
$\mu_0$ & 3.95$\times10^{-5}$ kg/m-s\\
$\boldsymbol{g}$ & (0,0,-9.8) m/$\mathrm{s}^2$\\
Injection pressure & 2.03$\times10^{9}$ Pa\\
Open boundary pressure & 3.075$\times10^{7}$ + 1.025$\times10^{4}$(depth) Pa\\
\hline
Nodes & 23126 \\
Elements & 43098 \\
\hline
\end{tabular}
\label{table:LeakyWellParams}
\end{table}

\subsection{Computational efficiency}
In this section, we evaluate both elements from the perspective of computation time and memory usage. The most significant drawback pointed out for the RT0 element is that it is not scalable as the problem size grows. This is due to the RT0 element formulation having variables that exist on the element edges, which in turn leads to more equations in the linear system. Simply put, a node unknown can be shared by many elements, where as an edge or face unknown can be shared by only two elements. The fewer elements that can share an unknown, the larger the system. Although the problems shown in this work are only two dimensional, effects of this drawback are present. Table \ref{table:Timing} shows the number of unknowns, assembly time, solve time, and memory usage of both elements for all of the examples above. Note that for each problem the number of unknowns is substantially higher for the RT0 element, roughly 30\%. If this were the only consideration, the VMS element would be the clear winner, but taking into account the computation time as well, the VMS element takes substantially longer to assemble and solve the linear system. This is due to the additional stabilization terms that must be computed at each step. As expected, the time savings decrease as the problem size increases since at some point the problem size offsets the faster computation for the RT0 element.

The same is true for the memory usage. For the two smaller problems, the memory usage is identical, but for the leaky well problem, the VMS element requires less memory. The results suggest that for small to moderately sized problems (less than $100,000$ degrees of freedom) the RT0 element is more efficient, but for large problems the VMS element is more efficient.

\begin{table}[htb!]
\caption{Timing and memory usage data}
\centering
\begin{tabular}{l l l l l l l}
\hline
 & \multicolumn{2}{l}{Multilayer}  & \multicolumn{2}{l}{Cylinder Inclusion} & \multicolumn{2}{l}{Leaky well} \\
Quantity & RT0 & VMS & RT0 & VMS & RT0 & VMS\\
\hline
Total unknowns &  2481 & 1764 &  2801 & 1956 &  132447 &  92504 \\
Assembly time/time step (s) & 0.0323 & 0.05162 & 0.0345 &  0.05746 & 2.943 & 4.246 \\
System solve/time step (s) & 0.0327 & 0.04056 & 0.0156 & 0.04959 & 2.564 &  2.320 \\
Required memory storage (MB) & 233 & 233 & 227 & 227 & 421 & 396 \\
\hline
Problem size savings &-- &28.9\% &-- &30.1\% & --& 30.2\%\\
Time savings & 29.5\% &-- &53.2\% &-- &16.1\% & --\\
Memory savings & --& --&-- &-- &-- &5.9\%\\
\hline
\end{tabular}
\label{table:Timing}
\end{table}